# Asymptotic sieve for primes

By John Friedlander and Henryk Iwaniec

## 1. Heuristics and statement of results

For a long time it was believed that sieve methods might be incapable of achieving the goal for which they had been created, the detection of prime numbers. Indeed, it was shown [S], [B] to be inevitable that, in the sieve's original framework, no such result was possible although one could come tantalizingly close. This general limitation is recognized as the "parity problem" of sieve theory.

More recently, beginning with the work [IJ], this goal has in certain cases become possible by adapting the sieve machinery to enable it to take advantage of the input of additional analytic data. There have been a number of recent developments in this regard, for example [H], [DFI], [FI], and several recent works of R.C. Baker and G. Harman. The story however is far from finished.

In this paper we consider a sequence of real nonnegative numbers

$$\mathcal{A} = (a_n) \tag{1.1}$$

with the purpose of showing an asymptotic formula for

$$S(x) = \sum_{n\leqslant x} a_n \Lambda(n) \tag{1.2}$$

where $\Lambda(n)$ denotes the von Mangoldt function. The sequence $\mathcal{A}$ can be quite thin although not arbitrarily so. Setting

$$A(x) = \sum_{n\leqslant x} a_n \tag{1.3}$$

we require that $A(x)$ be slightly larger than $A(\sqrt{x})$; precisely

$$A(x) \gg A(\sqrt{x})(\log x)^2 \ . \tag{1.4}$$

*JF was supported in part by NSERC grant A5123 and HI was supported in part by NSF grant DMS-9500797.

It is easy to predict the main term for $S(x)$ by using the formula

$$\Lambda(n) = -\sum_{d|n} \mu(d) \log d \ .$$

Inserting this in (1.2) we arrange $S(x)$ as

$$S(x) = \sum_{d\leqslant x} \lambda_d \, A_d(x)$$

where $\lambda_d = -\mu(d) \log d$ and $A_d(x)$ are the congruence sums

$$A_d(x) = \sum_{\substack{n\leqslant x \\ n\equiv 0 \pmod d}} a_n \ . \tag{1.5}$$

Hence the problem reduces to that of the asymptotic evaluation of the sums $A_d(x)$. In practice there is no problem to obtain a bound of the correct order of magnitude with considerable uniformity. We assume that

$$A_d(x) \ll d^{-1}\tau(d)^8 A(x) \tag{1.6}$$

uniformly in $d \leqslant x^{1/3}$. Of course this bound is crude and we need more precise formulas if we are insisting on asymptotic results. We assume that we can approximate $A_d(x)$ by writing

$$A_d(x) = g(d)A(x) + r_d(x) \ , \tag{1.7}$$

so $g$ may be regarded as a density function and $r_d(x)$ may be regarded as a rather small remainder term. In the leading term we assume that $g$ is multiplicative and satisfies

$$0 \leqslant g(p) < 1, \quad g(p) \ll p^{-1}, \tag{1.8}$$

for all primes. Moreover, some regularity in the distribution of $g(p)$ is required. We express this by assuming that for $y \geqslant 2$

$$\sum_{p\leqslant y} g(p) = \log\log y + c + O\left((\log y)^{-10}\right) \tag{1.9}$$

with some constant $c$.

Now we can write by (1.7)

$$S(x) = H(x)A(x) + R(x)$$

where

$$H(x) = -\sum_{d\leqslant x} \mu(d)g(d) \log d \tag{1.10}$$

and $R(x)$ is the remainder term

$$R(x) = \sum_{d\leqslant x} \lambda_d \, r_d(x) \ . \tag{1.11}$$

The hypothesis (1.9) allows us (see (2.4)) to extend (1.10) to the complete series getting by partial summation,

$$H(x) = H + O\left((\log x)^{-1}\right) \tag{1.12}$$

where $H$ is the positive constant

$$H = -\sum_{d} \mu(d)g(d)\log d\ . \tag{1.13}$$

That this is positive follows since the series is also given by the infinite product

$$H = \prod_{p} (1 - g(p)) \left(1 - \tfrac{1}{p}\right)^{-1}\ . \tag{1.14}$$

Hence we conclude that

$$S(x) = HA(x)\left\{1 + O\left((\log x)^{-1}\right)\right\} + R(x)\ . \tag{1.15}$$

This would finish the job if only one could ignore the remainder term $R(x)$, but in general this is nothing but a wishful dream. A primary reason is that it is impossible in practice to control the individual error terms $r_d(x)$ uniformly for $d \leqslant x$.

We assume that the error terms satisfy

$$\sum_{d \leqslant D} \mu^2(d)|r_d(t)| \leqslant A(x)(\log x)^{-2^{22}} \tag{R}$$

for all $t \leqslant x$ with some $D = D(x)$ such that

$$x^{\frac{2}{3}} < D(x) < x\ . \tag{R1}$$

The assumption (R) is quite realistic since $D(x)$ is not necessarily required to be near $x$. But such a choice of $D(x)$ necessitates some additional hypothesis by means of which one would be able to avoid an encounter with the extremely large moduli.

We assume an estimate for bilinear forms of the following type:

$$\sum_{m} \Big| \sum_{\substack{N<n\leqslant 2N\\ mn\leqslant x}} \gamma(n)\mu(mn)a_{mn} \Big| \leqslant A(x)(\log x)^{-2^{22}} \tag{B}$$

for every $N$ with

$$\Delta^{-1}\sqrt{D} < N < \delta^{-1}\sqrt{x} \tag{B1}$$

for some $\delta = \delta(x) \geqslant 2$ and $\Delta = \Delta(x) \geqslant 2$, and where the coefficients are given by

$$\gamma(n) = \gamma(n, C) = \sum_{d|n, d\leqslant C} \mu(d)\ . \tag{B2}$$

This is required for every $C$ with

$$1 \leqslant C \leqslant xD^{-1} . \tag{B3}$$

*Remarks.* Since $C$ is relatively small the coefficients $\gamma(n, C)$ behave more or less as constant; indeed we have $\gamma(n, 1) = 1$. A thorough examination of $\gamma(n, C)$ can be found in [DIT].

Finally, for simplicity we assume that the sequence $\mathcal{A}$ is supported on squarefree numbers, that is

$$a_n = 0 \qquad \text{if } \mu(n) = 0 . \tag{1.16}$$

In Section 9 we remove this assumption at the expense of strengthening slightly the hypothesis (R) and introducing a few extra minor conditions.

Our main result is

THEOREM 1. *Assuming the above hypotheses, we have*

$$\sum_{p \leqslant x} a_p \log p = HA(x) \left\{ 1 + O\left( \frac{\log \delta(x)}{\log \Delta(x)} \right) \right\} , \tag{1.17}$$

*where the implied constant depends only on the function $g$.*

In practice (B) can be established in the range (B1) for $\delta = (\log x)^\alpha$ and $\Delta = x^\eta$ with some positive constants $\alpha$, $\eta$. For these choices the error term in (1.17) becomes $O(\log \log x / \log x)$ and this cannot be improved by refining our arguments. If $\log \delta(x) \gg \log \Delta(x)$ then the error term exceeds the main term and (1.17) follows from the upper bound sieve. Therefore for the proof we may assume that $x \geqslant \Delta(x)^A$, $\Delta(x) \geqslant \delta(x)^A$, and $\delta(x) \geqslant A$ for any fixed positive constant $A$. Note also that it suffices to prove the Theorem with $\Delta(x)$ in the inequality (B1) replaced by $\Delta(x)^A$ (and similarly for $\delta(x)$).

Most of the assumptions made above are standard ones from sieve theory. We have tried to choose them so that they are not difficult to verify in practice and have made no attempt to choose a minimal set. For example, some of the information about the multiplicative function $g$ is stronger than necessary. The assumption (1.9) is as strong as the prime number theorem; it could be weakened considerably but it is not worth the effort. What needs to be said is that (1.9) is not crucial for breaking the parity problem.

In the classical sieve machinery all of the above assumptions but (B) are present in one form or another and, of these, (R) is the fundamental one. As already mentioned, with (R) alone one cannot capture the primes, even for $D(x) = x^{1-\varepsilon}$, which is the case treated by Bombieri [B]. It is clear that $D(x) = x$ cannot happen in practice. Actually, for thin sequences $\mathcal{A}$ one cannot expect (R) to hold with $D(x) > A(x)$ and for such sequences the classical sieve does correspondingly worse.

By introducing into the sieve axioms the additional hypothesis (B) we shall be able to resolve the "parity problem" and detect primes. Our theorem gives conditions under which one can guarantee the expected number of primes in the sequence $\mathcal{A}$. Note that in case $D(x) = x^{1-\varepsilon}$ very little is required in the range (B1) of the bilinear forms, namely

$$x^{\frac{1}{2}-\varepsilon'} < N < x^{\frac{1}{2}-\varepsilon}$$

for some $\varepsilon' > \varepsilon > 0$ with $\varepsilon' \to 0$ and $\varepsilon/\varepsilon' \to 0$, to produce the expected asymptotic formula. For thin sequences, $D(x)$ must be less than $A(x)$ and then the required range (B1) is somewhat larger but nevertheless we still retain the asymptotics for primes under reasonable assumptions. We could arrange our combinatorial arguments somewhat differently to come up with (B1) for smaller $N$ but it would be applicable only for sequences sufficiently dense.

Note that the restriction $N < \delta^{-1}\sqrt{x}$ in (B1) makes (B) realistic whereas requiring (B) for $N = \sqrt{x}$ would not be, for reasons similar to the prohibition of $D(x) = x$ in (R). Were we to make this seemingly slight extension in the range of (B), we would be able to give a simpler proof, but of a result with virtually no application. We overcome the additional difficulty occurring in the short range $\delta^{-1}\sqrt{x} \leqslant N \leqslant \sqrt{x}$ with the aid of a device used by Bombieri in [B]. Also, our stipulation of the lower bound restriction $N > \Delta^{-1}\sqrt{D}$ in (B1) is essential; indeed by narrowing this slightly to $N > \sqrt{D}$ we would not be able to break the parity problem.

It is worth mentioning that the source of cancellation in the bilinear form in (B) comes from the sign changes of the Möbius function $\mu(mn)$ in the inner sum so there is no need to express this sum in terms of the remainder $r_{mn}(x)$, as in some other works, cf. [I]. There is also the Möbius function $\mu(d)$ in the coefficient $\gamma(n, C)$. Here $d$ must be quite a bit smaller than $N$ to ensure that $\mu(d)$ does not completely neutralize $\mu(n)$. By (B1-B3) we know that $d < C < x\Delta N D^{-3/2}$, so our hypothesis (B) can be realistic only if $D$ is somewhat larger than $x^{2/3+\varepsilon}$.

To get an idea of the role of (B) in enabling us to detect primes let us recall a well-known example of Selberg which draws attention to the failure of the classical set-up. Let $\mathcal{A} = (a_n)$ with $a_n = \frac{1}{2}(1 + \lambda(n))$ where $\lambda$ is the Liouville function so $a_n$ is the characteristic function of the integers having an even number of prime factors, in particular $a_p = 0$ for every $p$. This sequence satisfies (R) with $D = x^{1-\varepsilon}$ but fails to satisfy (B). Even if we were to refine (B) replacing $a_{mn}$ by $r_{mn}(x)$ it would not be satisfied by $\mathcal{A}$. Indeed we have

$$2r_d(x) = \lambda(d) \sum_{\ell \leqslant \frac{x}{d}} \lambda(\ell) + O(1).$$

Thus $\mu(d)r_d(x)$ has constant sign in long intervals, namely for all $d$ satisfying $\frac{x}{L+1} < d \leqslant \frac{x}{L}$ with any integer $L < \sqrt{x}$.

Of course such a sequence is rather artificial and we certainly expect that the sieve given here is quite capable of settling, for instance, the twin prime problem. The stumbling block is that in this case we have no idea how to prove that the relevant sequence satifies the condition (B). Even proving (R) to such a high level is currently beyond reach. Fortunately there are other interesting sequences for which these hypotheses can be verified.

There are certainly other options for the introduction of bilinear form estimates to obtain primes with the sieve, see for example [DFI], but it is not yet clear what is the optimal general strategy if indeed there is one. Our choice was motivated by the desire to treat thin sequences, in particular $a_n$ being the number of representations $n = a^2 + b^4$. In this case (R) was established in [FI] with $D(x) = x^{3/4-\varepsilon}$. The bilinear form bound (B) valid in the range $x^{1/4+\varepsilon} < N < x^{1/2}(\log x)^{-A}$ is given in [FI2] and this more than suffices to prove the asymptotic formula

$$\sum_{a^2+b^4\leqslant x} \Lambda(a^2+b^4) \sim 4\pi^{-1}\kappa x^{\frac{3}{4}},$$

where $a$ and $b$ run over positive integers and $\kappa = \Gamma\left(\frac{1}{4}\right)^2/6\sqrt{2\pi}$. In Section 10, with such applications in mind, we give a modified theorem which incorporates a technical condition to make the axiom (B) easier to verify.

*Acknowledgment.* We are happy to thank the Institute for Advanced Study for pleasant working conditions and Enrico Bombieri for helpful conversations. We thank Andrew Granville, Mark Watkins, and the referee for pointing out inaccuracies in the previous draft and making suggestions which led us to improvements. J.F. also wishes to express thanks to Macquarie University for their hospitality during part of this work.

## 2. Technical reductions

What we actually need for the proof of Theorem 1 are the following hypotheses:

$$\text{(R}'\text{)}\qquad \sum_{d\leqslant D} \mu^2(d)\tau_5(d)|r_d(t)| \ll A(x)(\log x)^{-3}$$

$$\text{(B}'\text{)}\qquad \sum_m \tau_5(m) \Big| \sum_{\substack{N<n\leqslant 2N\\ mn\leqslant x}} \gamma(n)\mu(mn)a_{mn} \Big| \ll A(x)(\log x)^{-3}$$

in the same notation and ranges as in (R) and (B). Here the implied constant depends only on $g$.

In this section we derive (R′) and (B′) from (R) and (B), respectively, using the crude bound (1.6). As soon as (R′), (B′) are established the hypotheses (R), (B) and (1.6) are no longer required.

We begin with a lemma which is reminiscent of a result due to Wolke [W].

LEMMA 1. *Fix $k \geqslant 2$. Any squarefree integer $n \geqslant 1$ has a divisor $d \leqslant n^{1/k}$ such that*

$$\tau(n) \leqslant (2\tau(d))^k \ .$$

*Proof.* Write $n = p_1 \dots p_r$ with $p_1 < \cdots < p_r$ and put $d = p_1 \dots p_{[r/k]}$; thus $d \leqslant n^{1/k}$ and

$$2\tau(d) = 2^{1+[r/k]} \geqslant 2^{r/k} = \tau(n)^{1/k}$$

completing the proof. □

We use Lemma 1 for $k = 3$. First we estimate the following sums:

$$\sum_{d\leqslant x} \tau_5(d)^{\frac{4}{3}} A_d(x) \leqslant \sum_{d\leqslant x} \tau_{10}(d) A_d(x) = \sum_{n\leqslant x} \tau_{11}(n) a_n$$
$$\leqslant \sum_{n\leqslant x} \tau(n)^4 a_n \leqslant \sum_{d\leqslant x^{1/3}} (2\tau(d))^{12} A_d(x)$$

by the lemma. Then by (1.6) this is bounded by $O\left(A(x)V(x)\right)$ where

$$V(x) = \sum_{d\leqslant x}{}^{\flat} d^{-1}\tau(d)^{20} \ll (\log x)^{2^{20}}.$$

Here, and hereafter, $\sum^{\flat}$ indicates that the summation is restricted to squarefree integers. Therefore

$$\sum_{d\leqslant x} \tau_5(d)^{\frac{4}{3}} A_d(x) \ll A(x)(\log x)^{2^{20}} \ . \tag{2.1}$$

Similarly, but using (1.9),

$$\sum_{d\leqslant x}{}^{\flat} \tau_5(d)^{\frac{4}{3}} g(d) \leqslant \prod_{p\leqslant x} (1 + 11g(p)) \ll (\log x)^{11} \ . \tag{2.2}$$

Combining (2.1) and (2.2) we infer that

$$\sum_{d\leqslant x}{}^{\flat} \tau_5(d)^{\frac{4}{3}} |r_d(t)| \ll A(x)(\log x)^{2^{20}} \ . \tag{2.3}$$

Finally (R′) follows from (R) and (2.3) by Hölder's inequality.

The derivation of (B′) from (B) is similar.

We shall also need, and conclude this section by proving, the following consequence of (1.9):

$$\sum_{\substack{d\leqslant y\\(d,\nu)=1}}\mu(d)g(d)\ll\sigma_\nu(\log y)^{-6}\tag{2.4}$$

uniformly in $\nu\geqslant 1$, $y\geqslant 2$, where

$$\sigma_\nu=\prod_{p\mid\nu}\left(1+\tfrac{1}{\sqrt{p}}\right).\tag{2.5}$$

Let $2\leqslant z\leqslant\sqrt{y}$ and $P(z)$ be the product of the primes $p<z$. We split the sum (2.4), say $G_\nu(y)$, as follows:

$$\begin{aligned}G_\nu(y)&=\sum\sum_{\substack{mn\leqslant y,(mn,\nu)=1\\m\mid P(z),(n,P(z))=1}}\mu(m)\mu(n)g(m)g(n)\\&=\sum_{\substack{m\leqslant\sqrt{y}\\m\mid P(z),(m,\nu)=1}}\mu(m)g(m)\sum_{\substack{n\leqslant y/m\\(n,\nu P(z))=1}}\mu(n)g(n)\\&\quad+\sum_{\substack{n\leqslant\sqrt{y}\\(n,\nu P(z))=1}}\mu(n)g(n)\sum_{\substack{\sqrt{y}<m\leqslant y/n\\m\mid P(z),(m,\nu)=1}}\mu(m)g(m).\end{aligned}$$

Estimating the outer sums by

$$\sideset{}{^\flat}\sum_{n\leqslant\sqrt{y}}g(n)\leqslant\prod_{p\leqslant\sqrt{y}}(1+g(p))\ll\log y$$

we get

$$G_\nu(y)\ll(|G_\nu(w,z)|+H(y,z))\log y$$

for some $w$ with $\sqrt{y}\leqslant w\leqslant y$, where

$$G_\nu(w,z)=\sum_{\substack{n\leqslant w\\(n,\nu P(z))=1}}\mu(n)g(n)$$

and

$$H(y,z)=\sum_{\substack{m>\sqrt{y}\\m\mid P(z)}}g(m).$$

We estimate $H(y,z)$ by using Rankin's trick as follows:

$$\begin{aligned}H(y,z)&\leqslant y^{-\varepsilon}\sum_{m\mid P(z)}g(m)m^{2\varepsilon}=y^{-\varepsilon}\prod_{p<z}\left(1+g(p)p^{2\varepsilon}\right)\\&\leqslant y^{-\varepsilon}\prod_{p<z}(1+8g(p))\ll\exp\left(-\frac{\log y}{\log z}\right)(\log z)^8\end{aligned}$$

by choosing $\varepsilon = (\log z)^{-1}$. To estimate $G_\nu(w, z)$ we consider

$$F(w, z) = \sum_{\substack{n \leqslant w \\ (n, P(z))=1}} \mu(n) f(n)$$

where $f$ is another multiplicative function which satisfies the hypothesis (1.9). By (1.9) we deduce that

$$\left| \sum_{u<p\leqslant v} f(p) - \sum_{\substack{u<p\leqslant v \\ p\nmid \nu}} g(p) \right| \leqslant h(z)$$

for any $v \geqslant u \geqslant z$ with

$$h(z) \ll \sum_{p \geqslant z, p|\nu} p^{-1} + (\log z)^{-10} \ll \sigma_\nu (\log z)^{-10}.$$

We also have the inequality

$$|F(w, z) - G_\nu(w, z)| \leqslant h(z) \sum_{mn \leqslant w} f(m) g(n)$$

which can be seen from the identity

$$f(p_1 ... p_r) - g(p_1 ... p_r) = \sum_{1 \leqslant j \leqslant r} f(p_1 ... p_{j-1}) \left( f(p_j) - g(p_j) \right) g(p_{j+1} ... p_r).$$

Therefore

$$G_\nu(w, z) = F(w, z) + O\left( h(z) (\log w)^2 \right).$$

It remains to estimate $F(w, z)$. To this end we take the simple function $f(n) = n^{-1}$ for which we know that

$$F(w, z) \ll z^{-c} + \exp\left( -\frac{\log w}{\log z} \right) (\log z)^6.$$

This can be proven by a standard contour integration using the bound for the corresponding zeta function

$$\sum_{(n, P(z))=1} \mu(n) f(n) n^{-s} = \zeta(s+1)^{-1} \prod_{p<z} (1 - p^{-1-s})^{-1} \ll (\log(|s|+1))(\log z)^3$$

valid for

$$\mathrm{Re} s \geqslant -\min \left\{ \frac{c}{\log(|s|+2)}, \frac{1}{\log z} \right\},$$

where $c$ is a positive constant. Collecting the above results we arrive at

$$G_\nu(y) \ll \sigma_\nu (\log z)^{-10} (\log y)^3 + \exp\left( -\frac{\log y}{2 \log z} \right) (\log y)^9.$$

We choose $\log z = (\log y)(\log \log y)^{-2}$ obtaining (2.4).

## 3. Combinatorial identities

For $f : \mathbb{N} \to \mathbb{C}$, an arithmetic function, we split

$$f(n \leqslant z) = \begin{cases} f(n) & \text{if } n \leqslant z \\ 0 & \text{if } n > z \end{cases}$$

$$f(n > z) = \begin{cases} 0 & \text{if } n \leqslant z \\ f(n) & \text{if } n > z. \end{cases}$$

We have

$$f(n > z) = \sum_{bc|n} \mu(b) f(c > z).$$

Split the Möbius function $\mu(b) = \mu(b \leqslant y) + \mu(b > y)$ getting

$$f(n > z) = \sum_{bc|n} \mu(b \leqslant y) f(c > z) + \sum_{bc|n} \mu(b > y) f(c > z).$$

In the first sum insert $f(c > z) = f(c) - f(c \leqslant z)$ getting the identity (a variant of that of Vaughan),

$$f(n > z) = \sum_{b|n} \mu(b \leqslant y) F\left(\tfrac{n}{b}\right) - \sum_{bc|n} \mu(b \leqslant y) f(c \leqslant z) + \sum_{bc|n} \mu(b > y) f(c > z),$$

where $F = f * 1$. Split the last sum, say $f(n; y, z)$, into three subsums

$$f_1(n; y, z) = \sum_{bc|n} \mu(b > sy) f(c > sz),$$

$$f_2(n; y, z) = \sum_{bc|n} \mu(sy \geqslant b > y) f(c > z),$$

$$f_3(n; y, z) = \sum_{bc|n} \mu(b > sy) f(sz \geqslant c > z).$$

Here $s$ can be any number $> 1$. We choose $s$ to be the power of 2 in the interval

$$\frac{\Delta\sqrt{x}}{8\delta\sqrt{D}} < s \leqslant \frac{\Delta\sqrt{x}}{4\delta\sqrt{D}} .$$

Denote

$$F(n; y) = \sum_{b|n} \mu(b \leqslant y) F\left(\tfrac{n}{b}\right) ,$$

$$F(n; y, z) = \sum_{bc|n} \mu(b \leqslant y) f(c \leqslant z).$$

We have

$$f(n > z) = F(n; y) - F(n; y, z) + \sum_{j=1}^{3} f_j(n; y, z). \tag{3.1}$$

We apply the integral operator

$$I_h(Y, Z) = \int_Y^{eY} \int_Z^{eZ} h(y, z)(yz)^{-1} dy\, dz$$

to obtain an identity with smoothed parts. Note that if $h$ is independent of one or both variables the integration over the relevant variable(s) does nothing. When we apply this operator to (3.1) we denote the results by changing $y, z$ to $Y, Z$ as appropriate. For example for the truncated function $f(n > z)$ this operator gives

$$f(n > Z) = \begin{cases} f(n) & \text{if } n > eZ \\ f(n) \log \frac{n}{Z} & \text{if } Z < n \leqslant eZ \\ 0 & \text{if } n \leqslant Z. \end{cases}$$

We introduce this integral operator over the variables $y, z$ in order to facilitate their separation in various forthcoming multiple sums.

We shall apply the above considerations for $f = \Lambda$ in which case $F = L$ with $L(n) = \log n$. Thus we have from (3.1),

$$\Lambda(n > Z) = L(n; Y) - L(n; Y, Z) + \sum_{j=1}^{3} \Lambda_j(n; Y, Z). \tag{3.2}$$

We shall take $n \leqslant x$ and make the choices

$$Y = Z = \Delta^{-1}\sqrt{D}\ . \tag{3.3}$$

Note that when $y$ changes over the segment $Y < y \leqslant eY$ then $sy$ satisfies the inequalities $\frac{1}{8}\delta^{-1}\sqrt{x} < sy < \delta^{-1}\sqrt{x}$. The same holds for the variable $z$.

The obvious fact that, for $n$ squarefree, the left-hand side of (3.2) is supported on primes is lost when we look instead on the right-hand side (and this is the intention of the formula). In order to keep partial track of this information we introduce an upper-bound sieve $\{\lambda_\nu,\ \nu \leqslant \Delta\}$ of level $\Delta$ and having $|\lambda_\nu| \leqslant 1$. Such a sieve can be obtained, for example, by restricting $\mu(\nu)$ to certain numbers $\nu \leqslant \Delta$, including $\nu = 1$, having the property

$$\rho_n = \sum_{\nu | n} \lambda_\nu \geqslant 0\ .$$

Clearly $\Lambda(n > Z) = \rho_n \Lambda(n > Z)$, for squarefree $n$ since $Z > \Delta$. Therefore,

$$S(x, Z) = \sum_{n \leqslant x} a_n \Lambda(n > Z) = \sum_{n \leqslant x} a_n \rho_n \Lambda(n > Z) \tag{3.4}$$
$$= T(x; Y) - T(x; Y, Z) + \sum_{j=1}^{3} S_j(x; Y, Z)$$

where

$$T(x; Y) = \sum_{n \leqslant x} a_n \rho_n L(n; Y)$$
$$T(x; Y, Z) = \sum_{n \leqslant x} a_n \rho_n L(n; Y, Z)$$
$$S_j(x; Y, Z) = \sum_{n \leqslant x} a_n \rho_n \Lambda_j(n; Y, Z), \quad j = 1, 2, 3 \ .$$

The sum in (1.2)

$$S(x) = \sum_{p \leqslant x} a_p \log p$$

is closely approximated by $S(x; Z)$; precisely we have

$$0 \leqslant S(x) - S(x; Z) \leqslant \sum_{p \leqslant eZ} a_p \log p \ll A(x)(\log x)^{-1} \tag{3.5}$$

by (1.4). It remains to treat the various sums on the right side of (3.4).

## 4. Evaluation of $T(x; Y)$

We are able to treat the sum $T(x; y)$ for individual $x, y$. We split this by means of $\log \frac{n}{b} = \log n - \log b$, getting

$$T(x; y) = T_1(x; y) + T_2(x; y).$$

We have

$$T_1(x; y) = \sum_{n \leqslant x} a_n \rho_n \log n \sum_{\substack{b \mid n \\ b \leqslant y}} \mu(b) \tag{4.1}$$
$$= \sum_{b \leqslant y} \mu(b) \sum_{\substack{n \leqslant x \\ n \equiv 0 \pmod b}} a_n \rho_n \log n.$$

Here the inner sum is, by partial summation,

$$V_b(x) \log x - \int_1^x V_b(t) \frac{dt}{t} \tag{4.2}$$

where

$$V_b(x) = \sum_{\substack{n\leqslant x\\ n\equiv 0 \pmod b}} a_n \sum_{\nu|n} \lambda_\nu \tag{4.3}$$
$$= \sum_{\nu} \lambda_\nu A_{[\nu,b]}(x) = \sum_{\nu} \lambda_\nu \left\{ g([\nu,b])A(x) + r_{[\nu,b]}(x) \right\}.$$

We have $g([\nu,b]) = g(\nu)g\left(\frac{b}{(\nu,b)}\right)$ and

$$\sum_{b\leqslant y} \mu(b) g\left(\frac{b}{(\nu,b)}\right) = \sum_{d|\nu} \mu(d) \sum_{\substack{b\leqslant \frac{y}{d}\\ (b,\nu)=1}} \mu(b)g(b).$$

This last inner sum is, by (2.4), bounded by $O\left(\sigma_\nu (\log x)^{-6}\right)$. By this and (4.2) we deduce that the contribution to (4.1) coming from the main term in (4.3) is

$$\ll A(x)(\log x)^{-5} \sum_{\nu} |\lambda_\nu| g(\nu)\tau(\nu)\sigma_\nu \ll A(x)(\log x)^{-3}.$$

The remainder terms in (4.3) make a contribution to (4.1) which is

$$\ll (\log x) \sum_{\nu} \sum_{b} |\lambda_\nu \mu(b) r_{[\nu,b]}(x)|$$
$$\ll (\log x) \sideset{}{^\flat}\sum_{d<\Delta y} \tau_3(d)|r_d(x)| \ll A(x)(\log x)^{-2}$$

by (R$'$) since $\Delta y < D$. Thus we obtain the bound $T_1(x;y) \ll A(x)(\log x)^{-2}$.

We now consider

$$T_2(x;y) = -\sum_{b\leqslant y} \mu(b) \log b \sum_{\nu} \lambda_\nu A_{[\nu,b]}(x).$$

We insert the approximation (1.7) for $A_{[\nu,b]}(x)$. The remainder terms give a contribution which may be estimated in the same way as those for $T_1(x;y)$ giving the same bound $\ll A(x)(\log x)^{-2}$. Here, however, the main term provides also the main term for the total sum $S(x)$. This is

$$-A(x) \sum_{\nu} \lambda_\nu g(\nu) \sum_{b\leqslant y} \mu(b) g\left(\tfrac{b}{(\nu,b)}\right) \log b\ . \tag{4.4}$$

We extend the summation to all $b$ making an error $\ll A(x)(\log x)^{-3}$, by (2.4). Then, we evaluate the complete inner sum as

$$\sum_{\eta|\nu} \mu(\eta) \sum_{(b,\nu)=1} \mu(b)g(b)\log \eta b = \sum_{\eta|\nu} \mu(\eta) \sum_{(b,\nu)=1} \mu(b)g(b)\log b$$

and this is equal to $-H$ if $\nu = 1$ and equal to zero otherwise. Thus the sum (4.4), and also $T_2(x,y)$ are given by $HA(x) + O\left(A(x)(\log x)^{-2}\right)$.

Combining this with the bound for $T_1(x,y)$ we obtain

$$T(x,y) = HA(x) + O\left(A(x)(\log x)^{-2}\right). \tag{4.5}$$

## 5. Estimation of $T(x;Y,Z)$

It again suffices to estimate the contribution for given $y$, $z$ and here this is given by the sum

$$T(x;y,z) = \sum_{b\leqslant y}\mu(b)\sum_{c\leqslant z}\Lambda(c)\sum_{\nu}\lambda_\nu A_{[\nu,bc]}(x).$$

We insert the approximation (1.7) for $A_{[\nu,bc]}(x)$. For the contribution coming from the main terms, arguing as with $T_1(x;y)$, we are now led to estimate the sum

$$\sum_{b\leqslant y}\mu(b)\sideset{}{^\flat}\sum_{\substack{c\leqslant z\\ c\nmid b}}\Lambda(c)\,g\left(\tfrac{bc}{(\nu,bc)}\right) = \sum_{b\leqslant y}\mu(b)\,g\left(\tfrac{b}{(\nu,b)}\right)\sideset{}{^\flat}\sum_{\substack{c\leqslant z\\ c\nmid b}}\Lambda(c)\,g\left(\tfrac{c}{(\nu,c)}\right)$$

and by (2.4), (1.8) this is $\ll \tau(\nu)\sigma_\nu(\log y)^{-6}(\log \nu z)$. Summing over $\nu$, as with $T_1(x;y)$ we get a contribution $\ll A(x)(\log x)^{-3}$.

The contribution coming from the remainder terms is bounded by

$$\sum_\nu|\lambda_\nu|\sideset{}{^\flat}\sum_{b\leqslant y}\sideset{}{^\flat}\sum_{\substack{c\leqslant z\\ c\nmid b}}\Lambda(c)|r_{[\nu,bc]}(x)| \leqslant (\log yz)\sum_\nu\sideset{}{^\flat}\sum_{d\leqslant yz}\left|\, r_{[\nu,d]}(x)\,\right| \ll A(x)(\log x)^{-2}$$

by (R$'$) since $yz \leqslant e^2\Delta^{-2}D < \Delta^{-1}D$. From the two estimates we conclude that

$$T(x;y,z) \ll A(x)(\log x)^{-2}\ . \tag{5.1}$$

## 6. Estimation of $S_1(x;Y,Z)$

As before this may be estimated for individual $y$ and $z$. We have

$$S_1(x;y,z) = \sum\sum\sum_{\substack{bcd\leqslant x\\ b>sy,c>sz}}\mu(b)\Lambda(c)\rho_{bcd}a_{bcd}\ .$$

Note that the conditions of summation imply that $c$ is located in a short interval (in the logarithmic scale), namely

$$\delta^{-1}\sqrt{x} < c < \delta\sqrt{x}\ , \tag{6.1}$$

and $d$ is relatively small, namely

$$d < \delta^2 \ . \tag{6.2}$$

Therefore there is enough room for $b$ out of which to create a nice variable. By dropping $\mu(b)$ we deduce that

$$\begin{aligned} |S_1(x;y,z)| &\leqslant \sum_c \Lambda(c) \sum_d \sum_{n\leqslant x, cd|n} \rho_n a_n \\ &= \sum_c \Lambda(c) \sum_d \sum_\nu \lambda_\nu A_{[cd,\nu]}(x) \ . \end{aligned}$$

Of course, we have lost the possibility of detecting cancellation due to the sign changes of $\mu(b)$, but we have gained the congruence sums which can be evaluated asymptotically due to (R$'$). Note that $[cd,\nu] \leqslant cd\nu < \delta^3\Delta\sqrt{x} < D$. We insert (1.7) and estimate the resulting remainder by (R$'$) getting

$$|S_1(x;y,z)| \leqslant A(x)L(x)M(x) + O\left(A(x)(\log x)^{-2}\right) \tag{6.3}$$

where $L(x)$, $M(x)$ come from the main terms, namely

$$L(x) = \sum_c{}^{\flat} \Lambda(c) g(c),$$

$$M(x) = \sum_d{}^{\flat} \sum_\nu \lambda_\nu g([d,\nu]) \ .$$

Here $c, d$ run over the segments (6.1) and (6.2) respectively. By (1.8) we infer that

$$L(x) \ll \log \delta \ . \tag{6.4}$$

For the estimation of $M(x)$ we write $g([d,\nu]) = g(d)g(\nu/(\nu,d))$ and note that

$$\sum_\nu \lambda_\nu g(\nu/(\nu,d)) \geqslant 0 \ .$$

The latter is a property of any upper-bound sieve, and it holds for any multiplicative function $f(\nu)$ with $0 \leqslant f \leqslant 1$ in place of $g(\nu/(\nu,d))$. This positivity property allows us to introduce the factor $(\delta^2/d)^\varepsilon$ into $M(x)$ and extend the summation to all $d$ (Rankin's trick) getting

$$\begin{aligned} M(x) &\leqslant \delta^{2\varepsilon} \sum_\nu \lambda_\nu g(\nu) \sum_d{}^{\flat} g(d/(\nu,d)) d^{-\varepsilon} \\ &= \delta^{2\varepsilon} \prod_p (1 + g(p)p^{-\varepsilon}) \sum_\nu \lambda_\nu g(\nu) h(\nu) \end{aligned}$$

where $h(\nu)$ is the multiplicative function with

$$h(p) = (1+p^{-\varepsilon})(1+g(p)p^{-\varepsilon})^{-1} < 2(1+g(p))^{-1} \ .$$

Note that $g(p)h(p) < 1$, so the sieve theory applies giving

$$\sum_{\nu} \lambda_\nu g(\nu)h(\nu) \ll \prod_{p<\Delta} (1 - g(p)h(p)) = \prod_{p<\Delta} (1 - g(p))(1 + g(p)p^{-\varepsilon})^{-1} .$$

Hence

$$M(x) \ll \delta^{2\varepsilon} \prod_{p\geqslant\Delta} (1 + g(p)p^{-\varepsilon}) \prod_{p<\Delta} (1 - g(p)) .$$

We choose $\varepsilon = (\log \Delta)^{-1}$ getting by (1.9),

$$M(x) \ll \prod_{p<\Delta} (1 - g(p)) \ll (\log \Delta)^{-1} . \tag{6.5}$$

Combining (6.3), (6.4) and (6.5) we conclude that

$$S_1(x; y, z) \ll A(x) \frac{\log \delta}{\log \Delta} . \tag{6.6}$$

We remark that it is for this bound alone that we introduced the $\lambda_\nu$. Without doing so we would lose a logarithm, and this very tight bound cannot afford such a loss. The optimal choice of $\lambda_\nu$ belongs to the theory of the 2-dimensional sieve.

## 7. Estimation of $S_2(x; Y, Z)$

In this section and the next we finally encounter the sums where (due to contamination of the variables coming from the support of $\rho$) it is necessary to perform the integrated estimate. Whenever $a_n$ occurs we assume $n \leqslant x$ without writing this condition every time. We have

$$S_2(x; y, z) = \sum_{e} \sum_{y<b\leqslant sy} \sum_{c>z} \mu(b)\Lambda(c)\rho_{ebc}a_{ebc},$$

$$|S_2(x; y, z)| \leqslant (\log x) \sum_{k} \Big| \sum_{y<b\leqslant sy} \mu(b)\rho_{bk}a_{bk} \Big|$$

$$\leqslant (\log x) \sum_{\nu_1} \sum_{\nu_2} |\lambda_{\nu_1\nu_2}| \sum_{k} \Big| \sum_{\frac{y}{\nu_2}<b\leqslant\frac{sy}{\nu_2}} \mu(b)a_{\nu_1\nu_2 bk} \Big| .$$

In order to remove $\nu_2$ from the range of the inner summation we take advantage of the integration over $y$ (still not integrating over $z$).

$$S_2(x; Y, z) \leqslant (\log x) \sum_{\nu_1} \sum_{\nu_2} |\lambda_{\nu_1\nu_2}| \sum_{k} \int_{\frac{Y}{\Delta}}^{eY} \Big| \sum_{y<b\leqslant sy} \mu(b)a_{\nu_1\nu_2 bk} \Big| \frac{dy}{y} \tag{7.1}$$

$$\leqslant (\log x) \int_{\frac{Y}{\Delta}}^{eY} \sum_{\ell} \tau_3(\ell) \Big| \sum_{y<b\leqslant sy} \mu(b)a_{b\ell} \Big| \frac{dy}{y} .$$

By (B′) we conclude that

$$S_2(x;Y,z) \ll A(x)(\log x)^{-1} \ . \tag{7.2}$$

## 8. Estimation of $S_3(x;Y,Z)$

Recall that

$$S_3(x;y,z) = \sum_e \sum_{b>sy} \mu(b) \sum_{z<c\leqslant sz} \Lambda(c)\rho_{eb}a_{ebc} \ .$$

Here we have, for the first time, taken advantage of the fact that $\rho_{ebc} = \rho_{eb}$ since $c$ is prime. In essence this sum is of the same type as that in the previous section but with $b$ and $c$ interchanged. Thus by analogy to the previous argument, we put $c$ into the inner summation because its range is appropriate for (B′). However, the fact that $c$ is weighted by $\Lambda$ rather than $\mu$ necessitates splitting it further in order to obtain a factor with sign changes from which we may hope to detect cancellation. We write

$$\mu(c)\Lambda(c) = \sum_{m|c} \mu(m) \log \tfrac{m}{C} = \lambda^+(c) - \lambda^-(c)$$

(recall $c$ is squarefree) with $C = xD^{-1}$ where

$$\lambda^+(c) = \sum_{m|c} \mu(m) \log^+ \tfrac{m}{C}$$

$$\lambda^-(c) = \sum_{m|c} \mu(m) \log^+ \tfrac{C}{m} \ .$$

Note that the definition of $\lambda^-(c)$ is close to $\gamma(c,C)$ in (B2); precisely we have

$$\lambda^-(c) = \int_1^C \gamma(c,t)t^{-1}dt \ .$$

This gives a contribution to $S_3(x;y,z)$ which, by virtue of (B′) applied $\log_2 s$ times, satisfies

$$|S^-(x;y,z)| \leqslant \sum_\ell \tau_4(\ell) \Big| \sum_{z<c\leqslant sz} \mu(c)\lambda^-(c)a_{c\ell} \Big| \ll A(x)(\log x)^{-1} \ . \tag{8.1}$$

The contribution to $S_3(x;y,z)$ from $\lambda^+(c)$ is arranged as follows

$$\begin{aligned} S^+(x;y,z) &= \sum_e \sum_{b>sy} \mu(b)\rho_{eb} \sum\sum_{z<\ell m\leqslant sz} \mu(\ell)\big(\log^+ \tfrac{m}{C}\big)a_{eb\ell m} \\ &= \int_C^x \Bigg( \sum_{b>sy} \mu(b) \sum_{ebz<x} \rho_{eb} \sum_{\substack{\ell t<sz \\ eb\ell t<x}} \mu(\ell) \sum_{\substack{m>t \\ z<\ell m\leqslant sz}} a_{eb\ell m} \Bigg) \frac{dt}{t} \ . \end{aligned}$$

Put $d = eb\ell$, so $d < D$ and $d$ is squarefree because of (1.16). Here the inner sum is equal to $A_d(\min\{x, ebsz\}) - A_d(\max\{dt, ebz\})$. Applying (1.7) this inner sum becomes

$$g(d) \sum_{\substack{dt<n\leqslant x\\ ebz<n\leqslant ebsz}} a_n + r_d(\min\{x, ebsz\}) - r_d(\max\{dt, ebz\}) \ . \tag{8.2}$$

We first treat the contribution to $S^+(x; y, z)$ coming from the main term in (8.2). Recall that $d = eb\ell$. We get a significant cancellation from summation of $\mu(b)$. Note that $b$ ranges over the interval $b_{\min} < b < b_{\max}$, where $b_{\min} = \max\{sy, n/esz\}$ and $b_{\max} = \min\{n/ez, n/e\ell t\}$. This yields

$$\sum_{\substack{(b,e\ell)=1\\ b_{\min}<b<b_{\max}}} \mu(b)\rho_{eb} g(eb\ell) = g(e\ell) \sum_{\nu_1} \mu(\nu_1) g(\nu_1) \sum_{\nu_2|e} \lambda_{\nu_1\nu_2} \sum_{\substack{(b,e\ell\nu_1)=1\\ b_{\min}<b\nu_1<b_{\max}}} \mu(b) g(b)$$

and here, by (2.4), the inner sum is $O\left(\sigma_{e\ell\nu_1}(\log x)^{-6}\right)$. Summing over all the other variables trivially we find that the main term in (8.2) contributes to $S^+(x; y, z)$ an amount $S^*(x; y, z)$ which satisfies

$$S^*(x; y, z) \ll A(x)(\log x)^{-1} \ . \tag{8.3}$$

The remainder terms in (8.2) contribute an amount $S'(x; y, z)$ satisfying

$$|S'(x; y, z)| \leqslant \int_C^x \sideset{}{^\flat}\sum_{dt<x} \sum_{k|d, kz<x} \tau_4(k) \left| \, r_d\left(\min\{x, ksz\}\right) - r_d\left(\max\{dt, kz\}\right) \right| \frac{dt}{t} \ .$$

Change the variable of integration $t$ into $t/d$, integrate in $z$ over $Z < z < eZ$ and change $z$ into $z/k$ getting

$$|S'(x; y, Z)| \leqslant \int_Z^x \int_C^x \sideset{}{^\flat}\sum_{d<D} \tau_5(d) \left| \, r_d\left(\min\{x, sz\}\right) - r_d\left(\max\{t, z\}\right) \right| \frac{dt}{t}\frac{dz}{z} \ .$$

Now we can apply (R$'$) obtaining

$$S'(x; y, Z) \ll A(x)(\log x)^{-1} \ . \tag{8.4}$$

Collecting (8.1), (8.3) and (8.4) we conclude that

$$S_3(x; y, Z) \ll A(x)(\log x)^{-1} \ . \tag{8.5}$$

Combining (3.4), (3.5), (4.5), (5.1), (6.6), (7.2) and (8.5) we complete the proof of Theorem 1.

## 9. A generalization

Our condition that $a_n$ are supported on squarefree numbers may not be convenient in some applications. Therefore, in this section we establish a result which does not require this condition. Throughout we assume that $\mathcal{A} = (a_n)$ satisfies all the previous hypotheses except possibly for (1.16). We replace (1.16) by the following condition:

For all primes $p$ we have

$$0 \leqslant g(p^2) \leqslant g(p), \quad g(p^2) \ll p^{-2}. \tag{9.1}$$

We may remark that, although this condition was not needed before it would have been natural to include it in the form $g(p^2) = 0$ because $A_{p^2}(x) = 0$ due to (1.16).

Now we need to have some control on the size of the coefficients $a_n$. We assume that

$$\sum_{n \leqslant x} a_n^2 \leqslant x^{-\frac{2}{3}} A(x)^2 \ . \tag{9.2}$$

*Remarks.* The condition (9.2) holds for sequences $\mathcal{A} = (a_n)$ satisfying $a_n \leqslant n^{\varepsilon}$ and $A(x) \geqslant x^{\frac{2}{3}+2\varepsilon}$. Actually our argument requires a slightly weaker condition, for example the bound

$$\sum_{n \leqslant x} a_n^2 \leqslant x^{-1} A(x)^2 D(x)^{\frac{1}{2}} \tag{9.3}$$

would be sufficient.

We furthermore need a slightly stronger version of (R), namely for any $t \leqslant x$ we assume

$$\sideset{}{^3}\sum_{d < DL^2} |r_d(t)| \leqslant A(x) L^{-2} \tag{R$_3$}$$

where $L = (\log x)^{2^{24}}$ and $\sum^3$ restricts to cubefree moduli.

Theorem 2. *Assume the following hypotheses*: (1.4), (1.6), (1.8), (1.9), (B), (B1), (B2), (B3), (9.1), (9.2), (R$_3$) *and* (R1). *Then* (1.17) *holds true for* $\mathcal{A} = (a_n)$.

By analogy with the squarefree case, before proving Theorem 2 we show that (R$_3$) implies

$$\sideset{}{^3}\sum_{d < DL^2} \tau(d) |r_d(t)| \ll A(x) L^{-1} (\log x)^{2^{17}}. \tag{R$'_3$}$$

To this end we establish the following extension of Lemma 1.

Lemma 2. *Fix $k \geqslant 2$. Any $n \geqslant 1$ has a divisor $d \leqslant n^{1/k}$ such that*

$$\tau(n) \leqslant (2\tau(d))^{\frac{k \log k}{\log 2}} \ .$$

*Proof.* Write $n = \ell m$ where $\ell$ collects those primes occuring with exponent $\geqslant k$ and $m$ is made up of those having exponent $< k$. Thus $(\ell, m) = 1$. Put

$$d_\ell = \prod_{p^\alpha \| \ell} p^{[\alpha/k]}$$

so that $d_\ell \leqslant \ell^{1/k}$ and

$$\tau(d_\ell)^k = \prod_{p^\alpha \| \ell} \left(\left[\tfrac{\alpha}{k}\right] + 1\right)^k \geqslant \prod_{p^\alpha \| \ell} (\alpha + 1) = \tau(\ell) \ .$$

Write $m = p_1^{\alpha_1} \dots p_r^{\alpha_r}$ with $p_1 < \cdots < p_r$. Setting $d_m = p_1 \dots p_{[r/k]}$ we have $d_m < m^{1/k}$ and

$$(2\tau(d_m))^k \geqslant 2^r \geqslant \prod_j (\alpha_j + 1)^{\frac{\log 2}{\log k}} = \tau(m)^{\frac{\log 2}{\log k}} \ .$$

Take $d = d_\ell d_m$ completing the proof. ◻

We proceed to the proof of $(\mathrm{R}_3')$. By Lemma 2 with $k = 3$, the left-hand side is bounded by

$$27 \sum_{d \leqslant x^{\frac{1}{3}}}^{3} \tau(d)^5 |r_d(t)| \leqslant 27 \Bigg( \sum_{d \leqslant x^{\frac{1}{3}}}^{3} \tau(d)^{10} |r_d(t)| \Bigg)^{\frac{1}{2}} \Bigg( \sum_{d \leqslant x^{\frac{1}{3}}}^{3} |r_d(t)| \Bigg)^{\frac{1}{2}}$$

by Cauchy's inequality. In the first sum we replace $|r_d(t)|$ by $A_d(x) + g(d)A(x)$ and estimate as follows. By (1.6)

$$\sum_{d \leqslant x^{\frac{1}{3}}}^{3} \tau(d)^{10} A_d(x) \ll A(x) \sum_{d \leqslant x^{\frac{1}{3}}} d^{-1} \tau(d)^{18} \ll A(x) (\log x)^{2^{18}},$$

while, by (9.1) and (1.9),

$$\sum_{d \leqslant x^{\frac{1}{3}}}^{3} \tau(d)^{10} g(d) \leqslant \prod_{p \leqslant x} \Big( 1 + 2^{10} g(p) + 3^{10} g(p^2) \Big) \ll (\log x)^{2^{10}}.$$

For the second sum we apply $(\mathrm{R}_3)$. These give $(\mathrm{R}_3')$.

Now we are ready to prove Theorem 2. We shall apply Theorem 1 for the sequence $\tilde{\mathcal{A}} = (\tilde{a}_n)$ with

$$\tilde{a}_n = \mu^2(n) a_n \ . \tag{9.4}$$

Let $\tilde{g}(d)$ be the multiplicative function given by

$$\tilde{g}(p) = \frac{g(p) - g(p^2)}{1 - g(p^2)} \ . \tag{9.5}$$

Notice that $\tilde{g}(p)$ satisfies (1.8) and (1.9) due to (9.1). For any $d$ squarefree we write

$$\tilde{A}_d(x) = \sum_{\substack{n\leqslant x\\ n\equiv 0 \pmod d}} \tilde{a}_n = \tilde{g}(d)\tilde{A}(x) + \tilde{r}_d(x) \tag{9.6}$$

where $\tilde{A}(x) = \tilde{A}_1(x)$.

Although (B) for $\mathcal{A}$ and $\tilde{\mathcal{A}}$ coincide, it takes some effort to verify (R) for $\tilde{r}_d$, that is

$$R(t,D) = \sum_{d\leqslant D}^{\flat} |\tilde{r}_d(t)| \leqslant \tilde{A}(x)(\log x)^{-2^{22}} \ . \tag{9.7}$$

We shall derive this from ($R_3$). Detecting squarefree numbers by the Möbius formula $\mu^2(n) = \sum_{\nu^2|n}\mu(\nu)$ we get, for $d$ squarefree

$$\tilde{A}_d(x) = \sum_{\nu}\mu(\nu)A_{[\nu^2,d]}(x) = \tilde{g}(d)GA(x) + \sum_{\nu}\mu(\nu)r_{[\nu^2,d]}(x)$$

where

$$G = \prod_p \left(1-g(p^2)\right) \ , \tag{9.8}$$

so that

$$\tilde{g}(d)G = \sum_{\nu}\mu(\nu)g\left([\nu^2,d]\right) \ . \tag{9.9}$$

In particular,

$$\tilde{A}(x) = GA(x) + \sum_{\nu}\mu(\nu)r_{\nu^2}(x) \ . \tag{9.10}$$

Hence

$$\tilde{r}_d(x) = \sum_{\nu}\mu(\nu)\left(r_{[\nu^2,d]}(x) - \tilde{g}(d)r_{\nu^2}(x)\right) \ .$$

Therefore we have $R(t,D) \ll E\log x$, where

$$E = \sum_{\nu}^{\flat}\sum_{d\leqslant D}^{\flat}|r_{[\nu^2,d]}(t)| \ .$$

Note that every $k$ has at most $\tau(k)$ representations as $k = [\nu^2, d]$ and $k$ is cubefree if it has any. By ($R_3'$) it follows that the error terms $r_{[\nu^2,d]}(t)$ with $\nu \leqslant L$ and $d \leqslant D$ contribute

$$E_1 \ll A(x)L^{-1}(\log x)^{2^{17}} \ . \tag{9.11}$$

From the remaining terms we get

$$\begin{aligned} E_2 &\leqslant \sum_{\nu > L}^{\flat} \sum_{d \leqslant x}^{\flat} \left( A_{[\nu^2, d]}(x) + g\left([\nu^2, d]\right) A(x) \right) \\ &= E_{21} + E_{22} + O\left( A(x) (\log x) L^{-1} (\log L)^B \right) \end{aligned}$$

where $E_{21} + E_{22}$ is the splitting of the first term according to $\nu > L\sqrt{D}$ or $L < \nu \leqslant L\sqrt{D}$ and the second term is estimated by using (1.8), (1.9), and (9.1). Here $B$ depends only on $g$. We estimate as follows:

$$\begin{aligned} E_{21} &\leqslant \sum_{n \leqslant x} a_n \tau(n) \sum_{\substack{\nu^2 | n \\ \nu > L\sqrt{D}}} 1 \leqslant \Big( \sum_{n \leqslant x} a_n^2 \Big)^{\frac{1}{2}} \Big( \sum_{n \leqslant x} \tau(n)^3 \sum_{\substack{\nu^2 | n \\ \nu > L\sqrt{D}}} 1 \Big)^{\frac{1}{2}} \\ &\ll x^{-\frac{1}{3}} A(x) D^{-\frac{1}{4}} L^{-\frac{1}{2}} x^{\frac{1}{2}} (\log x)^{17} \leqslant A(x) L^{-\frac{1}{2}} (\log x)^{17} \end{aligned}$$

by (9.2) and (R1). For the second sum we apply Cauchy's inequality getting

$$E_{22}^2 \leqslant \Big( \sum_{L < \nu \leqslant L\sqrt{D}}^{\flat} A_{\nu^2}(x) \Big) \Big( \sum_{n \leqslant x} a_n \tau(n)^3 \Big) = S_1 S_2 \ ,$$

say. Here

$$S_1 = \sum_{L < \nu \leqslant L\sqrt{D}}^{\flat} \left( g(\nu^2) A(x) + r_{\nu^2}(x) \right) \ll A(x) L^{-1} (\log L)^B$$

by (1.8), (9.1), and ($R_3$). To estimate $S_2$ we apply Lemma 2 with $k = 3$ getting by (1.6)

$$\begin{aligned} S_2 = \sum_{n \leqslant x} a_n \tau(n)^3 &\leqslant \sum_{d \leqslant x^{1/3}} A_d(x) \left( 2\tau(d) \right)^{\frac{9 \log 3}{\log 2}} \\ &\ll A(x) \sum_{d \leqslant x} d^{-1} \tau(d)^{\frac{9 \log 3}{\log 2} + 8} \ll A(x) (\log x)^{3^9 + 8} \ . \end{aligned}$$

Gathering the above estimates we conclude the proof of (9.7), but with $A(x)$ in place of $\tilde{A}(x)$ (actually we saved at least a factor $\log x$ which is needed for clearing the implied constants below). Note, however, that along the above lines we have, by (9.10), also proved

(9.12) $$\tilde{A}(x) = G A(x) \left\{ 1 + O\left( (\log x)^{-1} \right) \right\} \ .$$

This formula gives (9.7) for the sequence $\tilde{\mathcal{A}}_d$ because $x$ is large and it gives (1.6) because

$$\tilde{A}_d(x) \leqslant A_d(x) \ll d^{-1} \tau(d)^8 A(x) \ll d^{-1} \tau(d)^8 \tilde{A}(x).$$

This completes the proof that $\tilde{\mathcal{A}}$ satisfies all the hypotheses of Theorem 1.

The constant (1.14) corresponding to $\tilde{g}$ may be expressed as

$$\tilde{H} = \prod_p (1 - \tilde{g}(p)) \left(1 - \tfrac{1}{p}\right)^{-1} = HG^{-1} \ , \tag{9.13}$$

with $G$ as in (9.8), so that (9.12) yields

$$\tilde{H}\, \tilde{A}(x) = HA(x) \left\{1 + O\left((\log x)^{-1}\right)\right\} \ . \tag{9.14}$$

Therefore the formula (1.17) for the sequence $\tilde{\mathcal{A}} = (\tilde{a}_n)$ coincides with that for $\mathcal{A} = (a_n)$, as it has to. This completes the proof of Theorem 2.

## 10. A reduction of the bilinear form

Our new axiom (B) in sieve theory can be made technically easier to verify in practice if the involved integer variables have no small prime divisors. Having this in mind, in this section we introduce an extra restriction into the inner sum in (B), namely

$$(n, \Pi) = 1 \tag{10.1}$$

where $\Pi$ is the product of all primes $p < P$ and $P$ is given subject to

$$2 \leqslant P \leqslant \Delta^{1/2^{35} \log\log x} \ . \tag{10.2}$$

Notice that this restriction is possible if $\Delta > (\log x)^{2^{35}}$. We assume in place of (B) that

$$\sum_m \Big| \sideset{}{^*}\sum_{\substack{N<n\leqslant 2N \\ mn \leqslant x}} \gamma(n)\mu(mn)a_{mn} \Big| \ll A(x)(\log x)^{-2^{26}} \ , \tag{B$^*$}$$

where $*$ restricts the inner summation to numbers $n$ free of prime factors $p < P$.

THEOREM 3. *Replace* (B) *by* (B$^*$) *in the assumptions of Theorem* 2. *Then* (1.17) *still holds.*

For the proof it is enough to show that (B$^*$) implies (B) with $\delta$, $\Delta$ replaced by $2\delta$, $2\Delta^2$ respectively. Furthermore we only need to establish (B) in the following integrated form

$$\int_N^{2N} \int_1^C \sum_m \Big| \sum_{w<n\leqslant 2w} \gamma(n,t)\mu(mn)a_{mn} \Big| \frac{dt}{t}\frac{dw}{w} \ll A(x)(\log x)^{-2^{22}} \tag{$\int$B}$$

for every $N$ with $\Delta^{-1}\sqrt{D} < N < \delta^{-1}\sqrt{x}$ in two cases $\gamma(n,t) = 1$ and

$$\gamma(n,t) = \sum_{d|n, d\leqslant t} \mu(d) \ . \tag{10.3}$$

To see this recall that (B) was used solely to estimate the sums $S_2$ and $S_3$ in Sections 7 and 8. In the case of $S_2$ we note that (7.1) contains such an integration with $y$ in place of $w$ and $\gamma(n,t)=1$. In the case of $S_3$ we note that (8.1) may be integrated over $z$ in place of $w$ and there $\lambda^-(c)$ is the integral of $\gamma(c,t)$ as required.

We write $n = n_0 n_1$ where $(n_0,\Pi)=1$ and $n_1 \mid \Pi$. The contribution to ($\int$B) from terms with $n_1 > \Delta$ is estimated trivially as

$$W_1 \leqslant (\log C) \sum\sum\sum_{\substack{mn_0n_1\leqslant x\\ n_1>\Delta}} \tau(n_0n_1)\mu^2(mn_0n_1)a_{mn_0n_1}.$$

Every such $n_1$ has a divisor $d$ with $\Delta < d \leqslant \Delta P$, so

$$W_1 \leqslant (\log C) \sum_{\substack{d\mid\Pi\\ \Delta<d\leqslant\Delta P}} \sideset{}{^\flat}\sum_{\substack{n\leqslant x\\ d\mid n}} a_n\tau_4(n)\ .$$

By Lemma 1 for $k=4$ we obtain $\tau_4(n) \leqslant \tau(n)^3 \leqslant (2\tau(\nu))^{24}$ for some $\nu|n$, $\nu \leqslant x^{\frac14}$ and hence by (1.6)

$$W_1 \leqslant (\log C) \sum_{\substack{d\mid\Pi\\ \Delta<d\leqslant\Delta P}} \sideset{}{^\flat}\sum_{\substack{\nu\leqslant x^{1/4}\\ (\nu,d)=1}} (2\tau(\nu))^{24} A_{\nu d}(x) \ll A(x)G_1G_2\log x$$

where $G_1$, $G_2$ are the following sums:

$$G_1 = \sideset{}{^\flat}\sum_{\nu\leqslant x} \nu^{-1}\tau(\nu)^{32} \ll (\log x)^{2^{32}},$$

$$G_2 = \sum_{\substack{d\mid\Pi\\ d>\Delta}} d^{-1}\tau(d)^8 \leqslant \Delta^{-\varepsilon}\sum_{d\mid\Pi} d^{\varepsilon-1}\tau(d)^8 = \Delta^{-\varepsilon}\prod_{p<P}\left(1+2^8p^{\varepsilon-1}\right)\ ,$$

by Rankin's trick. We choose $\varepsilon = (\log P)^{-1}$ so that $p^\varepsilon < 8p^{-\varepsilon}$ for every $p \leqslant P$ and obtain

$$\begin{aligned} G_2 &\leqslant \Delta^{-\varepsilon}\prod_p\left(1+2^{11}p^{-1-\varepsilon}\right) \leqslant \Delta^{-\varepsilon}\zeta(1+\varepsilon)^{2^{11}} \\ &\ll \Delta^{-\varepsilon}\varepsilon^{-2^{11}} = \exp\left(-\frac{\log\Delta}{\log P}\right)(\log P)^{2^{11}} \ll (\log x)^{-2^{34}}\ . \end{aligned}$$

Combining the above estimates we conclude that

$$W_1 \ll A(x)(\log x)^{-2^{33}}\ . \tag{10.4}$$

Next we estimate the contribution of terms with $n_1 \leqslant \Delta$. These give

$$W_0 \leqslant \int_N^{2N}\int_1^C \sum_m \sum_{\substack{n_1\mid\Pi\\ n_1\leqslant\Delta}} \Big| \sideset{}{^*}\sum_{w<n_0n_1\leqslant 2w} \gamma(n_0n_1,t)\mu(mn_0n_1)a_{mn_0n_1} \Big| \,\frac{dt}{t}\,\frac{dw}{w}\ .$$

In the case (10.3) we have

$$\gamma(n_0n_1, t) = \sum_{d_1|n_1} \mu(d_1)\gamma(n_0, t/d_1).$$

Insert this and change variables $t \to td_1$, $w \to wn_1$ getting

$$W_0 \leqslant \int_{N/\Delta}^{2N} \int_1^C \sum_m \tau_3(m) \mid \sum_{w<n\leqslant 2w}^{*} \gamma(n,t)\mu(mn)a_{mn} \mid \frac{dt}{t}\frac{dw}{w} .$$

The same bound holds in the case $\gamma(n,t) = 1$, in fact in stronger form since then the integration in $t$ is not needed and $\tau_3(m)$ can be replaced by $\tau(m)$. By the argument that gave (B) $\Rightarrow$ (B$'$) it follows that (B$^*$) implies

$$\sum_m \tau_3(m) \mid \sum_{N<n\leqslant 2N}^{*} \gamma(n)\mu(mn)a_{mn} \mid \ll A(x)(\log x)^{-2^{23}}$$

(the restriction $(n, \Pi) = 1$ only helps matters). Hence we conclude that

$$W_0 \ll A(x)(\log x)^{-2^{22}} . \tag{10.5}$$

Adding (10.4) to (10.5) we obtain ($\int$B) completing the proof of Theorem 3.

University of Toronto, Toronto, Canada
*E-mail address*: frdlndr@math.toronto.edu
Rutgers University, New Brunswick, NJ
*E-mail address*: iwaniec@math.rutgers.edu